\documentclass[oneside, 12pt]{amsart}

\usepackage{amscd,amssymb,enumerate, amsmath,graphicx}
\usepackage{mathrsfs}

\newcommand\mylabel[1]{\label{#1}}

\newtheorem{theorem}{Theorem}[section]

\newtheorem{lemma}[theorem]{Lemma}
\newtheorem{proposition}[theorem]{Proposition}
\newtheorem{corollary}[theorem]{Corollary}

\theoremstyle{definition}

\newtheorem*{acknowledgement}{Acknowledgement}

\theoremstyle{remark}

\DeclareFontFamily{U}{wncy}{}
\DeclareFontShape{U}{wncy}{m}{n}{<->wncyr10}{}
\DeclareSymbolFont{mcy}{U}{wncy}{m}{n}
\DeclareMathSymbol{\Sh}{\mathord}{mcy}{"58}
\newcommand{\ZZ}	{\mathbb{Z}}
\newcommand{\QQ}	{\mathbb{Q}}

\newcommand{\FF}	{\mathbb{F}}

\newcommand{\ideala}    {\mathfrak{a}}
\newcommand{\idealb}    {\mathfrak{b}}

\newcommand  {\shF}     {\mathscr{F}}

\newcommand  {\Ass}     {\operatorname{Ass}}

\newcommand  {\can}     {{\rm \text{can}}}

\newcommand  {\et}      {{\text{\rm et}}}
\newcommand  {\Et}      {{\text{\rm Et}}}

\newcommand  {\Frac}    {\operatorname{Frac}}

\newcommand  {\id}      {{\operatorname{id}}}

\newcommand  {\length}	{\operatorname{length}}

\newcommand  {\dirlim}  {\varinjlim}
\newcommand  {\invlim}  {\varprojlim}

\newcommand  {\lra}     {\longrightarrow}

\newcommand  {\maxid}   {\mathfrak{m}}

\newcommand  {\Nil}     {\operatorname{Nil}}

\newcommand  {\primid}  {\mathfrak{p}}

\newcommand  {\quadand} {\quad\text{and}\quad}

\newcommand  {\ra}      {\rightarrow}

\newcommand  {\rank}    {\operatorname{rank}}
\newcommand  {\red}     {{\operatorname{red}}}

\newcommand  {\Spec}    {\operatorname{Spec}}

\newcommand{\prc}	{\text{\rm{prc}}}

\def\mydate{\number\day\space\ifcase\month \or January\or February\or March\or 
April\or May\or June\or July\or
August\or September\or October\or November\or December\fi \space\number\year}

\DeclareFontFamily{U}{wncy}{}
\DeclareFontShape{U}{wncy}{m}{n}{<->wncyr10}{}
\DeclareSymbolFont{mcy}{U}{wncy}{m}{n}
\DeclareMathSymbol{\Sh}{\mathord}{mcy}{"58}
\begin{document}

\title[The $p$-radical closure]
      {The $p$-radical closure of  local noetherian rings}

\author[Stefan Schr\"oer]{Stefan Schr\"oer}
\address{Mathematisches Institut, Heinrich-Heine-Universit\"at,
40204 D\"usseldorf, Germany}
\curraddr{}
\email{schroeer@math.uni-duesseldorf.de}

%\keywords{}

\dedicatory{Final version, 5 May 2017}

\begin{abstract}
Given a local noetherian ring $R$ whose formal completion is integral,
we introduce and study the $p$-radical closure $R^\prc$.
Roughly speaking, this is the largest purely inseparable $R$-subalgebra
inside the formal completion $\hat{R}$.
It turns out that the finitely generated intermediate rings $R\subset A\subset R^\prc$
have rather peculiar properties. They can be used in  a systematic way to provide
examples of integral local rings 
whose normalization is non-finite,
that do not admit a resolution of singularities, and whose formal completion is
non-reduced.
\end{abstract}

\maketitle
\tableofcontents

%===========================================================
\section*{Introduction}
\mylabel{introduction}

In commutative algebra and algebraic geometry,
one frequently considers   \emph{integral closures} $R\subset R'$ of an integral noetherian domain $R$
with respect to    finite extensions $F\subset F'$ of the field of fractions $F=\Frac(R)$.
The whole theory of number fields and their rings of integers hinges on this process.
Moreover, the induced morphism $\Spec(R')\ra\Spec(R)$ can be regarded as a  higher-dimensional analog
of branched coverings of Riemann surfaces.

Indeed, under very general assumptions, the $R$-algebras $R'$ are finite, as in the case
of number rings and Riemann surfaces.
This holds, for example, if the   field extension $F\subset F'$ is separable, or if 
the ring $R$ is essentially of finite type over a field,
or a complete local noetherian ring.
This finiteness   is also a consequence of the defining conditions for \emph{excellent rings},
a class introduced by Grothendieck \cite{EGA IVb} that  is stable under
forming algebras of finite type, localizations and completions.
A  recent overview was given by Raynaud and Laszlo  in \cite{Illusie; Laszlo; Orgogozo 2014}, Expos\'e I.

However, the finiteness property does not hold for each and every noetherian ring, not even for all discrete valuation rings.
To my knowledge, the first counterexamples were  devised by
Akizuki \cite{Akizuki 1935} in characteristic zero, see also Reid's discussion   \cite{Reid 1995}, Section 9.5,
and Schmidt \cite{Schmidt 1936} in characteristic $p>0$. A counterexample
that is a discrete valuation ring
was given by Nagata (\cite{Nagata 1962}, Example (E3.3) on page 207):
the tensor product ring $R=k^p[[T]]\otimes_{k^p}k$, which is 
is a discrete valuation ring with $\hat{R}=k[[T]]$. Here 
$k$ is a field of positive characteristic $p>0$ with infinite $p$-degree, that is, $[k:k^p]=\infty$.
He also found intermediate discrete valuation rings $R\subset A\subset\hat{R}$
with    $\hat{A}=k[[T]]$
so that $A\subset\hat{A}$ is an infinite integral extension of finite degree $p$.
The study of intermediate rings $R\subset A\subset\hat{R}$ initiated by Nagata has a long tradition:
Bennett \cite{Bennett 1973} gave a detailed study,  via birational geometry,
of bad one-dimensional local noetherian rings. 
Olberding \cite{Olberding 2014} analyzed which discrete valuation rings appear 
as normalization of such rings.
Heinzer, Rotthaus and Wiegand studied various aspects of intermediate rings
$R\subset A\subset \hat{R}$   in a series of papers, among others in 
\cite{Heinzer; Rotthaus; Wiegand 1997}, \cite{Heinzer; Rotthaus; Wiegand 1999},
\cite{Heinzer; Rotthaus; Wiegand 2001}, \cite{Heinzer; Rotthaus; Wiegand 2015}.

Recall the following terminology:
An integral  domain $R$ is called \emph{japanese} if 
for every finite extension $F\subset F'$ of its field of fractions, the integral
closure $R'$ is a finite $R$-algebra. 
Now let us restrict our consideration to    local noetherian domains $R$.
Then $R'$ is a finite $R$-algebra provided the \emph{generic formal fiber} $\hat{R}\otimes_RF$ is reduced.
In fact,  Rees \cite{Rees 1951}  showed that the generic formal fiber is reduced if and only if 
all the intermediate rings $R\subset R'\subset\Frac(R)$ that are finitely generated $R$-algebras are japanese.
A comprehensive account on rings without finite normalization is given by  Olberding \cite{Olberding 2012}.
Recently, Koll\'ar \cite{Kollar 2015} reformulated the theory of normalizations in terms of pairs $(X,Z)$ 
consisting of noetherian scheme $X$ and a closed subscheme $Z\subset X$, which yields notions  that are preserved under formal completion
and work well even for schemes with non-reduced formal fibers.

The goal of this paper is to introduce and study the \emph{$p$-radical closure} $R\subset R^\prc$
for a local noetherian ring $R$ whose formal completion $\hat{R}$ is integral.
If $\hat{R}$ is normal, this equals the integral closure of $R$ with respect to the relative $p$-radical
closure of $\Frac(R)\subset\Frac(\hat{R})$, where $p\geq 1$ is the
characteristic exponent of the fields of fractions.
In some sense, it might be viewed as the  purely inseparable analogue of the henselization $R^h\subset\hat{R}$.
It would be   interesting to compute the $p$-radical closure in concrete  examples.
Here, however, our main goal is to uncover interesting general facts and formal consequences for rings having 
nontrivial   $R\subsetneqq R^\prc$.

It turns out that the \emph{intermediate rings} between $R$ and its $p$-radical closure $R^\prc$ have amazing properties,
and yield, in a rather systematic way,  examples of local noetherian rings with bad behavior.
Key features are as follows:
If $R\subsetneqq A\subset R^\prc$ is an intermediate ring that
is finite as  an  $R$-algebra, then its formal completion $\hat{A}$ contains nilpotent elements,
and  we actually have  $(\hat{A})_\red=\hat{R}$. The latter is our key observation, and most of the paper
hinges on this result. These rings $A$ are also  examples of 
local noetherian rings that do not admit resolutions of singularities.
In the one-dimensional situation, they must be non-normal, and
the  normalization map is non-finite. 

However, if $R\subsetneqq  B\subset R^\prc$ is an intermediate ring that is noetherian,
with reduced generic formal fibers, then $\hat{B}=\hat{R}$ holds and $R\subset B$ is not finite;
it follows  that the extension $R\subset B$ is faithfully flat,
and $B$ inherits regularity and Cohen--Macaulay properties from $R$.
It would be interesting to know under what circumstances the family of such
intermediate rings is cofinal. From this it would follow 
that the whole $p$-radical closure $R^\prc$ belongs to this family, at least if $R$ is regular.

In the case where $R$ is a discrete valuation ring, the theory simplifies a lot:
This is due to the Krull--Akizuki Theorem, 
and the fact that there is only one relevant formal fiber  $\hat{R}\otimes_R\kappa(\primid)$, 
namely the generic formal fiber.

\medskip
The paper is organized as follows:
In the first section, we review several well-known facts pertaining
to extensions and formal completions of local rings used throughout.
Since we have to cope with non-noetherian rings, we do not restrict our
discussion to the noetherian case.
In Section \ref{P-radical closure} we introduce the $p$-radical closure for local noetherian rings $R$
that are formally integral, and establish 
its basic properties.
The heart of the paper is Section \ref{intermediate rings}: Here we study the intermediate rings
between $R$ and its $p$-radical closure $R^\prc$.
In the final Section \ref{Discrete valuation rings}, we specialize our findings to the case of discrete valuation rings.

\begin{acknowledgement}
I wish to thank the referees  for several suggestions, which helped to improve the paper,
and for pointing out some bad noetherian local rings in \cite{Nagata 1962}.
\end{acknowledgement}

%===========================================================
\section{Extensions and completions of local rings}
\mylabel{extensions completions}

Here we recall some well-known basic facts on local rings and their formal completions
that will be used throughout the paper. Even though we are mainly interested in noetherian
rings, the non-noetherian rings are included in our discussion, because they may show up naturally
when it comes to taking integral closures.

A homomorphism $\varphi:R\ra A$ between local rings is  \emph{local}
if $\varphi^{-1}(\maxid_A)=\maxid_R$. 
An \emph{extension of local rings} is an injective local homomorphism $\varphi:R\ra A$ between local rings,
and we then usually write   $R\subset A$.
In general, one has the following two \emph{numerical invariants}
$$
e=\length_A(A/\maxid_RA)\quadand f=[\kappa(A):\kappa(R)],
$$
which we refer to as the \emph{ramification index} $e\geq 1$ and the \emph{residual degree} $f\geq 1$
for $\varphi:R\ra A$.
Here  $\kappa(R)=R/\maxid_R$ and $\kappa(A)=A/\maxid_A$ are the 
residue fields, and the numerical invariants are regarded as elements
from $\{1,2,\ldots\}\cup\{\infty\}$. 
If $R$ is integral, with field of fractions $F=\Frac(A)$,
there is another invariant  
$$
n=[A\otimes_RF:F]
$$
called the \emph{degree} $n\geq 0$ of the   homomorphism $R\ra A$. In case that the $R$-algebra $A$
flat, finitely presented and finite, that is, the underlying $R$-module is free of finite rank,
the famous formula $n=ef$ holds. 
These numerical invariants originate from the theory
of number fields and Riemann surfaces, and I find them quite useful in the general context.
Throughout, we are particularly interested in extensions of local rings with invariants
$e=f=1$. In  the context of valuation theory,
such extension are called \emph{immediate extensions}.

The \emph{formal completion} of the local ring $R$  is written as
$$
\hat{R}=\invlim_{n\geq 0} R/\maxid_R^n.
$$
This is another local ring, with maximal ideal $\maxid_{\hat{R}}=\maxid_R\hat{R}$
and residue field $\kappa(\hat{R})=\kappa(R)$ (\cite{AC 1-4}, Chapter III, \S2, No.\ 13, Proposition 19).
The canonical map $R\ra\hat{R}$, $a\mapsto (a,a,\ldots)$ is local,
with numerical invariants $e=f=1$. The local ring $R$ is called \emph{complete} if 
$R\ra\hat{R}$ is bijective.  
Clearly, the kernel of the canonical map is  the intersection
$\bigcap_{n\geq 0}\maxid_R^n$. It is trivial if $R$ is noetherian, or more generally if
 $R$, viewed as a topological ring with respect to the $\maxid_R$-adic topology,
 is   Hausdorff.
On the other hand, the kernel coincides with the maximal ideal if $R$ is 
an fppf-local ring, as studied  by Gabber and Kelly  in \cite{Gabber; Kelly 2015}, Section 3 and  myself in
\cite{Schroeer 2016}, Section 4.
Using \cite{AC 1-4}, Chapter III, \S2, No.\ 10, Corollary 5, we immediately get:

\begin{proposition}
\mylabel{formal completion noetherian}
The complete local ring $\hat{R}$ is noetherian if and only if the cotangent space
$\maxid_R/\maxid_R^2$ has finite dimension as vector space over the residue field $\kappa(R)=R/\maxid_R$.
\end{proposition}

\medskip
If $\varphi:R\ra A$ is a local homomorphism between local rings,
we get an induced local map $\hat{\varphi}:\hat{R}\ra\hat{A}$,
making the diagram
$$
\begin{CD}
\hat{R}	@>\hat{\varphi} >>	\hat{A}\\
@A\can AA		@AA\can A\\
R	@>>\varphi >	A
\end{CD}
$$
commutative. Consequently, the residual degree $f\geq 1$ for $R\ra A$ and $\hat{R}\ra\hat{A}$ coincide.
Furthermore:

\begin{proposition}
\mylabel{ramification index coincide}
If the fiber ring $A/\maxid_RA$ is noetherian,   
then the ramification indices $e\geq 1$
for the local maps $R\ra A$ and $\hat{R}\ra \hat{A}$ coincide.
\end{proposition}

\proof
The ramification indices in question are the lengths of the rings $B=A/\maxid_RA$
and $B'=\hat{A}/\maxid_{\hat{R}}\hat{A}=\hat{A}/\maxid_R\hat{A}$. The latter is complete, so
the canonical map $B\ra B'$ factors over
$\hat{B}$, and the resulting map $\hat{B}\ra B'$ is bijective.
Replacing $A$ by $B$, we thus reduce to the situation where $A$ is noetherian.

Now consider the complements $U\subset\Spec(A)$ and $U'\subset\Spec(\hat{A})$ of the closed point.
The resulting morphism $f:U'\ra U$ is quasicompact and faithfully flat. Let $M$ be a finitely generated $A$-module,
and $\shF$ be the quasicoherent sheaf on $U$ obtained by restricting the quasicoherent sheaf  
$\tilde{M}$ on $\Spec(A)$.
By fpqc-descent, the preimage $f^*(\shF)$ vanishes if and only if $\shF$ vanishes 
(\cite{SGA 1}, Expose VIII, Section 1, Corollary 1.3).
If follows that $M$ has finite length if and only if $\hat{M}=M\otimes_A\hat{A}$ has finite length.
In this situation,  $\maxid_A^nM=0$ for some $n\geq 0$, whence $\hat{M}=M$, and
one easily infers that $\length_A(M)=\length_{\hat{A}}(\hat{M})$.
Applying this for $M=A/\maxid_RA$ yields the assertion.
\qed

\medskip
The schematic fibers of the morphism $\Spec(\hat{R})\ra\Spec(R)$
are called the \emph{formal fibers} of  $R$.
The corresponding rings
$\hat{R}\otimes_R \kappa(\primid) $, where  $\primid\subset R$ are the prime ideals,
are called the \emph{formal fiber rings}. Note that these are noetherian rings, provided that $\hat{R}$ is noetherian.
They are endowed with the structure of an  algebra   over the  fields $\kappa(\primid)= R_\primid/\primid R_\primid$,
but these algebras are usually  not finitely generated.
The closed fiber equals the spectrum of the residue field $k=\hat{R}/\maxid_{\hat{R}}=R/\maxid_R$,
whence is of little interest.
Rather  important are the \emph{generic formal fibers}, which are given by the rings
$\hat{R}\otimes_R \kappa(\primid) $, where  $\primid\subset R$ are the minimal prime ideals.
Let us recall:

\begin{proposition}
\mylabel{reduced, irreducible, integral}
Let $R$ be a local noetherian ring. 
Then the  complete local scheme $\Spec(\hat{R})$ is
reduced or irreducible if and only if the respective property holds for the
local scheme $\Spec(R)$ and its generic formal fiber.
\end{proposition}

\proof
Suppose $\hat{R}$ is reduced. Then so is the subring $R\subset\hat{R}$.
Let $\primid\subset R$ be a minimal prime ideal, and write $S= R\smallsetminus\primid$ for 
the ensuing multiplicative system.  Being a localization of a reduced ring,
the generic formal $\hat{R}\otimes_R \kappa(\primid) = S^{-1}\hat{R}$ is reduced.

Conversely, suppose that  $R$ and its generic formal fiber is reduced.
We have to show that $\hat{R}$ contains
no embedded prime, and is regular in codimension zero.
According to \cite{Matsumura 1980}, Corollary in (9.B), each associated point
$\eta\in\Spec(\hat{R})$ lies over a generic point of $\Spec(R)$, and  
is generic in its formal fiber. Thus $\eta\in\Spec(\hat{R})$ is generic.
Moreover, the local ring corresponding to $\eta\in\Spec(\hat{R})$ is regular,
according to \cite{EGA IVb}, Corollary 6.5.2. Thus $\hat{R}$ is reduced.
The arguments for irreducibility are analogous, and left to the reader.
\qed

\medskip
Note that it may happen that $R$ is integral, yet there are embedded primes $\primid\subset\hat{R}$,
as examples of Ferrand and Raynaud reveal \cite{Ferrand; Raynaud 1970}.
Such behavior is ubiquitous, according to Lech \cite{Lech 1986}.
In this situation, however, we must have $\primid\cap R=0$, that is,
$\Ass(\hat{R})$ maps to the generic point of $\Spec(R)$.

We say that $R$ is \emph{formally reduced} if the equivalent conditions for reducedness of
the preceding Proposition hold. In \cite{Nagata 1962}, the term \emph{analytically unramified} was used.
We prefer the former locution, because ``unramifed'' now is used 
in algebraic geometry in a different way (\cite{EGA IVd}, \S 17).
We say that $R$ is \emph{formally integral} if the equivalent conditions for
integrality hold.
In \cite{Nagata 1962}, such rings are called \emph{analytically irreducible}.
The following observation, which is a consequence of  \cite{Nagata 1962}, (18.3) and (18.4),
applies in particular to the formal completion $R\subset\hat{R}$ 
of   formally integral local noetherian rings:

\begin{lemma}
\mylabel{fractions and completions}
Let $B\subset C$ be a faithfully flat extension of   domains.
Then we have $B=C\cap \Frac(B)$ as subsets of $\Frac(C)$.
\end{lemma}

If $R\subset A$ is an extension of local noetherian rings, the induced map $\hat{R}\ra\hat{A}$ stays
injective provided that the $R$-algebra $A$ is finite, by \cite{AC 1-4}, Chapter III, \S3, No.\ 4, Theorem 3. 
In general, however, the map is not injective,
for example if  $A$ but not $R$ is formally reduced.
We shall encounter such behavior later. 
The following is a situation in which injectivity is preserved without any finiteness assumptions.
The second assertion can be seen as a generalization of 
\cite{Olberding 2012}, Lemma 3.1, which dealt with discrete valuation rings:

\begin{proposition}
\mylabel{intermediate not finite}
Let $R\subset A$  be an extension of   local noetherian rings
with  numerical invariants $e=f=1$. Assume that $\dim(A)=\dim(R)$, and that $R$ is formally integral.
Then the induced map $\hat{R}\ra\hat{A}$
is bijective, and the $R$-algebra $A$ is faithfully flat.
Moreover, for any intermediate ring $R\subset R'\subsetneqq A$, the $R'$-algebra $A$ is not finite.
\end{proposition}

\proof
The local map $\hat{R}\ra\hat{A}$ has invariants $e=f=1$, by Proposition \ref{ramification index coincide}.
Hence the map is surjective, according to \cite{AC 1-4}, Chapter III, \S2, No.\ 9, Proposition 11.
Write   $\hat{A}=\hat{R}/\ideala$ for some ideal $\ideala\subset\hat{R}$.
Using
$$
\dim(A)=\dim(R)=\dim(\hat{R})\geq \dim(\hat{A})=\dim(A),
$$
we infer that $\dim(\hat{R})=\dim(\hat{A})$.
Since $\hat{R}$ is integral, we can apply Krull's Principal Ideal Theorem
and conclude that $\ideala=0$, whence the local map $\hat{R}\ra\hat{A}$ is bijective.
This ensures, by  \cite{AC 1-4}, Chapter III, \S3, No.\ 5, Proposition 10, that the $R$-algebra
$A$ is faithfully flat.

Now assume that for some intermediate ring $R\subset R'\subset A$ the $R'$-algebra $A$ is finite.
We have to verify that $R'=A$.
The map $\Spec(A)\ra\Spec(R')$ is   closed since it is finite. It is also dominant, whence surjective,
because the homomorphism $R'\ra A$ is injective.  Consequently, $R'$ is local,
and $R'\subset A$ and whence $R\subset R'$ are extensions of local rings.
By the Eakin--Nagata Theorem (\cite{Eakin 1968} or \cite{Nagata 1968}), the local ring $R'$ is noetherian.
The extension of local rings $R'\subset A$ has invariants $e=f=1$, because this holds for $R\subset A$.
Taking formal completions thus gives a finite extension of complete local rings 
$\hat{R}'\subset\hat{A}=A\otimes_{R'}\hat{R}'$ with
invariants $e=f=1$, whence the inclusion is an equality.
By faithfully flat descent, the inclusion $R'\subset A$ is an equality (\cite{SGA 1}, Expos\'e  VIII,
Corollary 1.3).
\qed

\medskip
Now suppose that our local noetherian ring $R$ is integral, with field of fractions
$F=\Frac(R)$. Then $R$ is called \emph{japanese} if for every finite extension 
$F\subset F'$, the integral closure $R\subset R'$ inside $F'$ is a finite $R$-algebra.
This automatically holds if the ring $R$ is normal and the field $F$ has characteristic zero 
(see \cite{Matsumura 1980}, Proposition 31.B).

A ring $R$ is called \emph{universally japanese} if each integral $R$-algebra of finite type
is japanese. Such rings are also known as  \emph{pseudo-geometric}, \emph{$J$-rings} or \emph{Nagata rings}.
A local ring $R$ is called \emph{excellent} if it is noetherian,
universally catenary, and has geometrically regular formal fibers;
such rings are universally japanese (see \cite{EGA IVb}, Scholie 7.8.3 (i) and (vi)).

%===========================================================
\section{The $p$-radical closure}
\mylabel{P-radical closure}

Throughout this section, $R$ denotes a local noetherian ring that is formally integral.
The extension of local rings  $R\subset\hat{R}$ induces an extension
of fields of fractions $\Frac(R)\subset\Frac(\hat{R})$.
We write $p\geq 1$ for the \emph{characteristic exponent} of these fields.
Recall that $p=1$ if the prime field is $\QQ$, and equals the characteristic $p\geq 2$ otherwise.
In the latter case, the canonical map $\ZZ\ra R\subset\Frac(R)$ yields an inclusion $\FF_p\subset R$.
Consider the subset 
$$
R^\prc=\{b\in\hat{R}\mid \text{$b^{p^\nu}\in R$ for some exponent $\nu\geq 0$}\}.
$$   
Clearly, this subset is a subring $R^\prc\subset\hat{R}$ containing $R$;  
we call it the \emph{$p$-radical closure} of $R$. Note that if $p=1$,
that is, the fields of fractions have characteristic zero, we have $ R^\prc=R$,
and the situation is of little interest. 
In general, set $F=\Frac(R)$ and let $F\subset F_\infty\subset \Frac(\hat{R})$ be the relative $p$-radical closure
in the sense of field theory.

\begin{proposition}
\mylabel{intersection and closure}
We have $R^\prc=\hat{R}\cap F_\infty$ inside the field of fractions $\Frac(\hat{R})$.
This is also the integral closure of $R$ with respect to $R\subset F_\infty$, provided
that $\hat{R}$ is normal.
\end{proposition}

\proof
The inclusion $R^\prc\subset \hat{R}\cap F_\infty$ is obvious.
For the reverse inclusion, suppose    $b\in\hat{R}$   is  a root 
of a polynomial $T^{p^\nu}-a$ with $a\in F$.
The task is to verify $a\in R$. Indeed, we have $a=b^{p^\nu}\in \hat{R}$,
and Lemma \ref{fractions and completions} tells us that $a\in \hat{R}\cap F=R$.
Now suppose that $\hat{R}$ is normal, and let  $R\subset R'\subset F_\infty$ be the integral closure.
Then $R'\subset\hat{R}$, because the latter is normal, and consequently
$R^\prc\subset R'\subset \hat{R}\cap F_\infty$. Since the right-hand side is contained
in $R^\prc$, the assertion follows.
\qed

\medskip
The $p$-radical closure is closely related to the generic formal fiber:

\begin{proposition}
\mylabel{geometrically reduced}
Consider the following three conditions:
\begin{enumerate}
\item The inclusion  $R\subset R^\prc$ is an equality.
\item The generic formal fiber of $R$  is geometrically reduced.
\item All finitely generated $R$-subalgebras $R'\subset\Frac(R)$  are japanese.
\end{enumerate}
Then  $ (ii) \Leftrightarrow (iii)$ holds. If $R$ is normal, we also have $(i) \Leftarrow (ii)$.
All three conditions are equivalent provided that $\hat{R}$ is normal and
the field extension $F\subset\Frac(\hat{R})$
can be written as  an purely inseparable extension followed by a separable extension.
\end{proposition}

\proof
The equivalence (ii) $\Leftrightarrow$ (iii) is a Theorem of Rees \cite{Rees 1951}.
Now suppose that $R$ is normal, and that (ii) holds.
Set $F=\Frac(R)$, and write $B=\hat{R}\otimes_RF$ for the generic formal fiber ring.
Suppose that $R^\prc\neq R$. Then there is some $b\in R^\prc$
with $b\not\in R$ but $a=b^p\in R$. Since $R$ is normal, we also have $b\not\in F$,
and   get a field extension $F\subsetneq F(a^{1/p})$. Hence the element $c=b\otimes 1-1\otimes a^{1/p}$
from the ring $B\otimes_F F(a^{1/p})$ is nonzero and   satisfies $c^p=0$,
thus $B$ is not geometrically reduced. 

Finally,  suppose that $\hat{R}$ is normal,   
that  the field extension $F\subset\Frac(\hat{R})$ can be written as  a purely inseparable extension $F\subset F'$ followed
by a separable extension $F'\subset\Frac(\hat{R})$, and that condition (i) holds, that is,  $R=R^\prc$. We shall verify (ii).
Clearly, the intermediate field $F'$ coincides with the the relative $p$-radical closure $F_\infty$.
First, we check that $F=F_\infty$: Suppose  $b/b'\in \Frac(\hat{R})$ is an element with 
$(b/b')^p=a/a'$ for some  $a,a'\in R$.
Then $(a'b/b')^p=aa'^{p-1}\in R\subset \hat{R}$.
Since $\hat{R}$ is normal, we must have $a'b/b'\in \hat{R}$, and thus $a'b/b'\in R^\prc=R$.
Consequently $b/b'=(a'b/b') /a'\in \Frac(R)=F$. This shows $F=F_\infty$,
and it follows that the field extension $F\subset \Frac(\hat{R})$ is separable.

Being a localization of the integral ring $\hat{R}$, the formal fiber ring
$B=\hat{R}\otimes_RF$ is integral, with $\Frac(B)=\Frac(\hat{R})$.
For any field extension $F\subset F'$, the inclusion $B\subset\Frac(\hat{R})$
induces an inclusion $B\otimes_FF'\subset\Frac(\hat{R})\otimes_FF'$.
The right-hand side is reduced, because the field extension $F\subset\Frac(\hat{R})$ is separable,
and thus the left-hand side is reduced as well.
Hence $B$ is geometrically reduced.
\qed

\medskip
Note that there are non-separable field extension $F\subset E$ whose relative $p$-radical closure
is $F_\infty=F$. There are examples where $F\subset E$ is finite (the \emph{exceptional  field extensions} from
\cite{A 4-7}, Chapter V,  Exercises 1 and 2 for \S 7), or where $F\subset E$ is relatively algebraically closed
(related to \emph{quasifibrations} $X\ra B$ of proper normal schemes with geometrically non-reduced generic fiber, 
compare the discussion in \cite{Schroeer 2010}).

The following functoriality property is immediate from the definition of the $p$-radical closure:

\begin{proposition}
\mylabel{functorial}
Let $R\subset A$ be an extension of local noetherian rings, both of which are formally integral.
Then the induced map $\hat{R}\ra\hat{A}$ sends $R^\prc$ to $A^\prc$.
\end{proposition}

\medskip
Note that the induced map $\hat{R}\ra\hat{A}$ is not necessarily injective;
this phenomenon was analyzed by H\"ubl \cite{Huebl 2001}.
If the map is injective,
one   may view both $R^\prc$ and $A$ as subrings inside $\hat{A}$.
Combining with Proposition \ref{geometrically reduced}, we get the following
property of the $p$-radical closure:
 
\begin{corollary}
\mylabel{minimal}
Assumptions as in the Proposition. Additionally suppose that the induced map $\hat{R}\ra\hat{A}$ is injective.
If the generic formal fiber of $A$
is geometrically reduced, then  we have $R^\prc\subset A$ as subrings inside $\hat{A}$. 
\end{corollary}

\medskip
Let us call a polynomial over a field \emph{purely inseparable} if it has exactly one
root in the algebraic closure.
With an additional hypothesis, we may regard the $p$-radical closure as
a \emph{purely inseparable algebraic closure}:

\begin{proposition}
Suppose the extension of local rings  $R^\prc\subset\hat{R}$ is flat.
Write $F=\Frac(R)$.
Then $R^\prc\subset\hat{R}$ is the set of all elements $b\in\hat{R}$ that satisfy an algebraic equation
$f(b)=0$ for some non-zero polynomial $f\in R[T]$ that is   purely inseparable as polynomial over $F$.
\end{proposition}

\proof
Clearly, every $b\in R^\prc$ satisfies the condition. Conversely, suppose that $b\in \hat{R}$
is the root of such a polynomial $f\in R[T]$.
Since this polynomial is purely inseparable 
it must be of the form $f(T)=c(T-b)^n$ for some non-zero $c\in R$ and  some integer $n\geq 1$.
Write $n=p^\nu m$ for some exponent $\nu\geq 0$ and $p\nmid m$. Then 
$$
f(T) = c(T^{p^\nu}-b^{p^\nu})^m = cT^{p^\nu m} + cmb^{p^\nu}T^{p^\nu(m-1)} + \ldots + cb^{p^\nu m}.
$$
Comparing coefficients and using  $m\in R^\times$, we see $cb^{p^\nu}\in R$.
The element $cb\in \hat{R}$ thus has $(cb)^{p^\nu}\in R$, hence $cb\in R^\prc$,
so $b=cb/c\in\hat{R}\cap\Frac(R^\prc)$.
Now we use that assumption that
$R^\prc\subset \hat{R}$ is flat. It follows from Proposition \ref{universal homeomorphism} below that it is 
faithfully flat; whence Lemma \ref{fractions and completions} ensures that $\hat{R}\cap\Frac(R^\prc)=R^\prc$. 
\qed

\medskip
Note that this interpretation makes no reference to the characteristic exponent;
perhaps this would lead to  a meaningful definition of $p$-radical closure for  arbitrary local noetherian rings,
which need not be formally integral.

Moreover, it reveals a striking    analogy between 
$p$-radical closure and  \emph{henselization}:
If $A$ is any local noetherian ring, with henselization $A^h$, the formal completion $\hat{A}$ is also henselian,
and the universal property of henselization gives inclusions $A\subset A^h\subset \hat{A}$,
see \cite{EGA IVd}, Theorem 18.6.6.
Note that if $A$ is formally irreducible, then $A^h$ is irreducible, which in turn means that
$A$ is unibranch.

If $A$ is formally normal and universally japanese, then the henzelization $A^h$ coincides
with the \emph{algebraic closure} of $A$ inside the formal completion $\hat{A}$,
according to \cite{Nagata 1962}, Theorem 44.1.
In other words, $A^h$ is the set of elements $b\in\hat{A}$ that satisfy an algebraic
equation $f(b)=0$ for some non-zero $f\in A[T]$; note that here $f$ is not necessarily monic.
This even holds if $A$ is merely integral, and its formal fibers are reduced,
and the generic formal fiber is normal, according to \cite{Kurke; Pfister; Roczen 1975}, Proposition 2.10.2.
This was further generalized to Hensel couples in \cite{Greco; Sankaran 1976}, Theorem 3,
see also Remark 2.  

We have the following variant, which says that for our formally integral local noetherian ring $R$ 
the henselization  is the
\emph{separable algebraic closure}   inside the  formal completion, under a  rather mild
assumption:

\begin{theorem}
\mylabel{henselization}
Suppose the henselization $R^h$ has 
geometrically connected   generic formal fiber. 
Write $F=\Frac(R)$.
Then $R^h\subset\hat{R}$ is the set of all elements $b\in\hat{R}$ that satisfy an algebraic equation
$f(b)=0$ for some non-zero polynomial $f\in R[T]$ that is  separable as polynomial over $F$.
\end{theorem}

\proof
First, we check that each $b\in R^h$ satisfies the conditions.
By definition of the henselization in \cite{EGA IVd}, Section 18.6 we have $b\in B_\primid$, where
$B$ is an \'etale $R$-algebra, and $\primid\subset B$ is a prime ideal lying over the maximal ideal $\maxid_R\subset R$.
After localization, we may assume that $B$ is integral, because $R$ is unibranch.
Let $g\in F[T]$ be the minimal polynomial for the element $b\in \Frac(B)$. Then $g$ is separable,
and multiplying with a suitable non-zero element $r\in R$ yields a polynomial $f(T)=rg(T)$ with coefficients in $R$.
The canonical inclusions $B_\primid\subset R^h\subset\hat{R}$ reveal that $b$ satisfies the desired conditions.
Note that this holds without the assumption on the generic formal fiber of $R^h$.

Conversely, suppose that we have $b\in \hat{R}$ satisfying $f(b)=0$ for some polynomial 
$$
f(T)=cT^m+ \lambda_{m-1}T^{m-1}+ \ldots+\lambda_0
$$
with coefficients in $R$ and leading coefficient $c\neq 0$, such that $f(T)$ is separable over $F$. 
Our goal is to show $b\in R^h$. 
Let us first assume that the algebraic equation is integral,
that is, $c=1$. 
Consider the   $R^h$-subalgebra $B=R^h[b]$ inside $\hat{R}$. The ring $B$ is integral, and
the $R^h$-algebra $B$ is finite.
Let $n=[B:R^h]$. It suffices to verify $n=1$, because then  
$b\in \hat{R}\cap \Frac(R^h)=R^h$, by Lemma \ref{fractions and completions}.
Indeed, it follows from  \cite{EGA IVd}, Theorem 18.6.6 that  the canonical map $\widehat{R^h}\ra \hat{R}$ is
bijective, so we may regard the inclusion $R^h\subset\hat{R}$ as the formal completion. 
Seeking a contradiction, we now suppose $n\geq 2$.

By assumption,  the generic formal fiber of $R^h$ is
geometrically connected.
It follows that $B\subset B\otimes_{R^h}\hat{R}$, $x\mapsto x\otimes 1$ corresponds to a  morphism on schemes
whose generic fiber is connected.
On the other hand, using that $B\subset\hat{R}$, we get a bijection
$$
(B\otimes_{R^h}B)\otimes_B\hat{R}\lra B\otimes_{R^h}\hat{R},\quad
x\otimes y\otimes z\longmapsto xy\otimes z.
$$
Since $\Frac(R^h)\subset\Frac(B)$ is a separable extension of degree $n\geq 2$,
it follows from the Galois Correspondence
that the inclusion   $B\subset B\otimes_{R^h}B$, $y\mapsto 1\otimes y$ corresponds to  a morphism of schemes
whose generic  fiber is disconnected.  
In turn,   the composite map 
$$
B\lra B\otimes_{R^h}B\lra (B\otimes_{R^h}B)\otimes_B\hat{R} \lra B\otimes_{R^h}\hat{R},
$$
which is nothing but the canonical inclusion $B\subset B\otimes_{R^h}\hat{R}$,
induces a morphism of schemes with disconnected generic fiber, contradiction.

We finally treat the general case, where the leading coefficient $c\in R$
an arbitrary non-zero element. Clearly, the element $b'=cb$ from $\hat{R}$ is a root for the monic polynomial
$T^m + c^0\lambda_{m-1}T^{m-1}+\ldots+c^{m-1}\lambda_0=c^{m-1}f(T/c)$, which is separable over $F$.
By the preceding paragraph, we have $b'\in R^h$, and in turn $b=b'/c\in\hat{R}\cap\Frac(R^h)$.
Since $R^h\subset\hat{R}$ is faithfully flat,   Lemma \ref{fractions and completions} applies, which  
gives $\hat{R}\cap\Frac(R^h)=R^h$. 
\qed

\medskip
Note that the   formal fibers of $R^h$ indeed are geometrically connected provided they are geometrically normal,
according to \cite{EGA IVd}, Theorem 18.9.1. I do not know an example where the generic formal fiber
of $R^h$ is geometrically disconnected.  This does not happen if $R^h$ has the approximation
property, and the latter implies, according to \cite{Rotthaus 1990}, that $R^h$ and whence $R$ are excellent.

Let me also make a comment on \'etale cohomology:
Given a scheme $X$, we denote by $\Et(X)$
the site of all \'etale morphisms $U\ra X$, endowed with the Grothendieck topology
whose covering families $(U_\alpha\ra U)_{\alpha\in I}$ are those with $\bigcup_{\alpha\in I}U_\alpha\ra U$ is surjective.
Let $X_\et$ be the resulting topos of sheaves.
For any abelian \'etale  sheaf $F$ on $X$, that is, an abelian sheaf on the site $\Et(X)$,
in other words  an object in the category $X_\et$, one gets the cohomology groups $H^r(X_\et,F)$.
Now set $X=\Spec(R)$ and $X^\prc=\Spec(R^\prc)$.
According to \cite{SGA 1}, Expos\'e IX, Theorem 4.10, together with footnote (5), the pullback functor
$\Et(X)\ra\Et(X^\prc)$, $U\mapsto U\times_XX^\prc$ is an equivalence of categories.
It induces a continuous map of topoi 
$$
\epsilon=(\epsilon_*,\epsilon^*): (X^\prc)_\et\lra X_\et,
$$
where the adjoint functors $\epsilon_*,\epsilon^*$ are equivalences 
of categories. In particular, computing \'etale cohomology on $X$ amounts to  the same as
computing it on $X^\prc$.

%===========================================================
\section{Intermediate rings}
\mylabel{intermediate rings}

Let $R$ be a local noetherian ring that is formally integral.
We now want to study intermediate rings between $R$ and its $p$-radical closure $R^\prc$.
Let us first record:

\begin{proposition}
\mylabel{universal homeomorphism}
For every intermediate ring $R\subset A\subset R^\prc$, the
induced morphisms 
$\Spec(R^\prc)\ra\Spec(A)$ and $\Spec(A)\ra\Spec(R)$ are  integral universal homeo\-morphisms.
The  extensions of local rings $R\subset A\subset R^\prc$ have numerical invariant $f=1$.
Moreover, the local rings $A$ are separated with respect to the $\maxid_A$-adic topology.
\end{proposition}

\proof
By definition, the morphisms  are integral. By the Going-Up Theorem, they are universally closed. 
Since the ring homomorphisms
are injective, the induced maps on affine schemes are dominant, whence surjective.
For each point $x\in\Spec(R)$, the fibers of $\Spec(R^\prc)\ra\Spec(R)$ and $\Spec(A)\ra\Spec(R)$
over $x$ are affine, and given by    algebras $C$ over $K=\kappa(x)$
where each $c\in C$ has some $c^{p^\nu}\in K$. If follows that $C\otimes_KK^{1/p^{\infty}}$
is a ring containing precisely one prime ideal, with trivial residue field extension. 
In light of \cite{EGA I}, Proposition 3.7.1, the morphisms
$\Spec(R^\prc)\ra\Spec(A)$ and $\Spec(A)\ra\Spec(R)$ are universally bijective. 
Being universally closed and universally bijective, they are universal homeo\-morphisms.

The assertion on residue fields holds, because all rings are contained in $\hat{A}$,
and the extension of local rings $A\subset\hat{A}$ has trivial residue field extension.
Finally, suppose we have some $f\in\bigcap_{n\geq 0}\maxid_A^n$. Using the
extension of local rings $A\subset \hat{R}$, we see that $f\in\bigcap_{n\geq 0}\maxid_{\hat{R}}^n=0$,
thus $A$ is separated.
\qed

\medskip
In particular, such $A$ are integral local rings,   the topological space $\Spec(A)$ is noetherian,
with $\dim(A)=\dim(R)$. Moreover,
$R\subset A$ is an extension of local rings, and we thus get an induced
map $\hat{R}\ra\hat{A}$ on formal completions. Note, however,
that there is no a priori reason that $A$ should be noetherian.

Now consider the direct system $A_\lambda\subset R^\prc$, $\lambda\in L$ of
all  finite $R$-subalgebras, that is, the underlying $R$-module is finitely generated.
The rings  $A_\lambda$ are noetherian by Hilbert's Basis Theorem.
Clearly, the direct system is filtered, 
and we have 
$$
R^\prc=\bigcup_{\lambda\in L} A_\lambda=\dirlim_{\lambda\in L}A_\lambda.
$$
The index set, viewed as an ordered set, has a smallest element $\lambda_{\min}\in L$.
To simplify notation, we denote this smallest element as $\lambda_{\min}=0$,
such that $A_0=R$.
Let us write $\hat{A}_\lambda=\widehat{A_\lambda}$ for the formal completions
of the local noetherian rings $A_\lambda$, and 
consider the \emph{reduced formal completions}
$$
(\hat{A}_\lambda)_\red = \hat{A}_\lambda/\Nil(\hat{A}_\lambda)
$$
and the composite map $\hat{R}\ra\hat{A}_\lambda\ra(\hat{A}_\lambda)_\red$.
The following is our key result:

\begin{theorem}
\mylabel{reduction}
For each $\lambda\in L$, the composite map $\hat{R}\ra(\hat{A}_\lambda)_\red$ is bijective.
\end{theorem}

\proof
The assertion is trivial if $R^\prc=R$. In light of Proposition \ref{geometrically reduced},
it thus suffices to treat the case that $F=\Frac(R)$ has positive characteristic $p>0$,
such that $R$ contains the prime field $\FF_p$.
First note that since $R\subset A_\lambda$ is finite, the canonical map
$$
A_\lambda\otimes_R \hat{R}\lra \hat{A}_\lambda,\quad x\otimes y\longmapsto xy
$$
is bijective (\cite{Matsumura 1980}, Theorem 55 on page 170). Now we use that the
$R$-algebra $A_\lambda$ is contained in $\hat{R}$, and observe that the resulting map
$$
(A_\lambda\otimes_R A_\lambda)\otimes_{A_\lambda} \hat{R}\lra A_\lambda\otimes_R \hat{R},\quad
(x\otimes y)\otimes z\longmapsto xy\otimes z
$$
is bijective. Clearly, the composite surjection
$$
(A_\lambda\otimes_R A_\lambda)\otimes_{A_\lambda} \hat{R}\lra 
A_\lambda\otimes_R \hat{R}\lra \hat{A}_\lambda\lra (\hat{A}_\lambda)_\red
$$
factors over $\tilde{R}=(A_\lambda\otimes_R A_\lambda)_\red\otimes_{A_\lambda}\hat{R}$.
% We thus get a commutative diagram
% $$
% \begin{CD}
% \tilde{R}	@>>> 	(\hat{A}_\lambda)_\red\\
% @AAA			@AAA\\
% \hat{R} 	@>>\id>	\hat{R}
% \end{CD}
% $$
Applying Lemma \ref{tensor product} below with $B=R$ and $C=A_\lambda$, we infer
that the mapping $A_\lambda\ra (A_\lambda\otimes_R A_\lambda)_\red$, $a\mapsto a\otimes 1$
is bijective. Consequently, we get a surjection
$$
\hat{R}\lra A_\lambda\otimes_{A_\lambda}\hat{R}\lra 
(A_\lambda\otimes_R A_\lambda)_\red\otimes_{A_\lambda}\hat{R}\lra (\hat{A}_\lambda)_\red.
$$
This map is given by  $a\mapsto 1\otimes a\mapsto (1\otimes 1)\otimes a\mapsto a$,
whence coincides with the the canonical map $\hat{R}\ra(\hat{A}_\lambda)_\red$.
Therefore, the latter is surjective.

Let $d\geq 0$ be the Krull dimension of the local noetherian ring $R$. 
In light of Proposition \ref{universal homeomorphism}, this is also
the dimension of $A_\lambda$, and thus of $(\hat{A}_\lambda)_\red$.
By assumption, the ring $\hat{R}$ is integral, whence
the surjective map  $\hat{R}\ra (\hat{A}_\lambda)_\red$ has trivial kernel,
by Krull's Principal Ideal Theorem.
\qed

\begin{corollary}
\mylabel{no resolution}
For each $\lambda\neq 0$, the  intermediate rings  $A_\lambda$   is  an integral  local
noetherian ring whose generic formal fiber  is non-reduced, and whose spectrum $\Spec(A_\lambda)$
does not admit a resolution of singularities.  
\end{corollary}

\proof
Suppose the generic formal fiber of  $A_\lambda$ is reduced. 
By Proposition \ref{reduced, irreducible, integral}, this means that $\hat{A}_\lambda$ is reduced.
In light of the Theorem, the map $\hat{A}\ra\hat{A}_\lambda$ is bijective.
By faithfully flat descent (\cite{SGA 1}, Expos\'e  VIII,
Corollary 1.3),  the inclusion $A\subset A_\lambda$ is an equality.
But this means $\lambda=0$, contradiction.

Thus the generic formal fiber of $A_\lambda$ is non-reduced, and in particular non-regular.
According to \cite{EGA IVb}, Proposition 7.9.3 this implies that the scheme
$\Spec(A_\lambda)$ admits no resolution of singularities.
\qed

\medskip
Using Proposition \ref{geometrically reduced}, we see that   normal local noetherian rings  $R$ 
whose generic formal fibers are integral but not geometrically reduced 
give  rise, in a rather systematic way, 
to  families of integral local noetherian rings 
without resolutions of singularities. In dimension one, this means:

\begin{corollary}
\mylabel{infinite normalization}
Suppose $\dim(R)=1$. Then for each $\lambda\neq 0$, the intermediate rings $A_\lambda$ 
is an integral    local noetherian ring with $\dim(A_\lambda)=1$
whose normalization $A_\lambda\subset A_\lambda'$ is not 
a finite extension.
\end{corollary}

\proof
If the normalization $A_\lambda\subset A_\lambda'$ is finite, then the normal one-dimensional ring $A'_\lambda$ is
noetherian, thus a Dedekind domain.
Hence  $\Spec(A'_\lambda)\ra\Spec(A_\lambda)$ is a resolution of
singularities, in contradiction to Corollary \ref{no resolution}.
\qed 

\medskip
Now let $R\subset B\subset R^\prc$ be an arbitrary intermediate ring.
If $B$ is not  a finite $R$-algebra, it may or may not be noetherian. 
Indeed, Nagata \cite{Nagata 1962} describes on page 207 examples constructed from   regular local rings of the form
$k^p[[x_1,\ldots,x_n]]\otimes_{k^p}k$, which yield non-noetherian $B$:
In Example 4 the ring  $R$ is integral of dimension two, and in Example 5 we have $R$ regular
of dimension three.
In any case,   $R\subset B$ is
an extension of local rings, hence we get an induced map on formal completions
$\hat{R}\ra\hat{B}$. Here $\hat{B}=\invlim B/\maxid_B^n$.

\begin{corollary}
\mylabel{image canonical map}
The image of the canonical map $B\ra(\hat{B})_\red$ is contained in the image
of $\hat{R}\ra(\hat{B})_\red$.
\end{corollary}

\proof
Let $b\in B$, and write $\bar{b}\in (\hat{B})_\red$ for its image.
Choose some index $\lambda\in L$ with $b\in A_\lambda\subset B$.
Consider the commutative diagram
$$
\begin{CD}
\hat{R}	@>>>	(\hat{A}_\lambda)_\red	@>>> 	(\hat{B})_\red\\
@AAA		@AAA				@AAA\\
R	@>>>	A_\lambda		@>>>	B.
\end{CD}
$$
According to the theorem, the upper left vertical map is surjective.
By construction, the class $\bar{b}$ lies in the image of the
upper right vertical map. It follows that $\bar{b}$ is also in the image of 
$\hat{R}\ra(\hat{B})_\red$.
\qed

\begin{corollary}
The reduced local ring $(\hat{B})_\red$ is noetherian. The local ring $\hat{B}$ is noetherian
if and only if its nilradical is finitely generated.
\end{corollary}

\proof
Since the local ring $\hat{B}$ is complete, so is the residue class ring $C=(\hat{B})_\red$.
It suffices to check that the latter has a finite cotangent space, according to Proposition 
\ref{formal completion noetherian}.
Since each vector from $\maxid_{C}/\maxid^2_{C}$ comes from some element in $B$,
Corollary \ref{image canonical map} reveals that it also comes from some element of $R$.
This implies that the canonical map $\maxid_{\hat{R}}/\maxid_{\hat{R}}^2\ra\maxid_{C}/\maxid^2_{C}$
is surjective. These vector spaces are finite-dimensional, because $\hat{R}$ is noetherian.

The condition in the second assertion is trivially necessary.
It is also sufficient: Suppose the nilradical $N\subset C$ is finitely generated.
Tensoring the exact sequence
$N\ra \maxid_{\hat{B}}\ra\maxid_C\ra 0$ of $\hat{B}$-modules with the residue field $k$ gives
an exact sequence
$$
N\otimes_{\hat{B}} k\lra \maxid_{\hat{B}}/\maxid^2_{\hat{B}}\lra \maxid_{C}/\maxid^2_{C}\lra 0.
$$
Consequently, the $k$-vector space in the middle is finite-dimensional.
Using Proposition \ref{formal completion noetherian}, we deduce that the complete local ring
$\hat{B}$ is noetherian.
\qed

\medskip
Our second main result deals with intermediate rings  that
are well-behaved:

\begin{theorem}
\mylabel{bijective}
Suppose  the intermediate ring $R\subset B\subset R^\prc$ is noetherian, and that its generic formal fiber is reduced.
Then the extension of local rings $R\subset B$ has invariants $e=f=1$, is faithfully flat,
and the induced map on formal completions $\hat{R}\ra\hat{B}$ is bijective.
Furthermore, if $B\neq R$ then the $R$-algebra $B$ is not finite.
\end{theorem}

\proof
We have  $f=1$  by Proposition \ref{universal homeomorphism}. In order to see $e=1$, we
first check that the induced map on cotangent spaces
\begin{equation}
\label{cotangent map}
\maxid_{R}/\maxid_{R}^2 \lra \maxid_{\hat{B}}/\maxid_{\hat{B}}^2
\end{equation}
is surjective. Let $\bar{b}\in \maxid_{\hat{B}}/\maxid_{\hat{B}}^2$. Since $B$
is noetherian, the map $\maxid_B/\maxid_B^2\ra \maxid_{\hat{B}}/\maxid_{\hat{B}}^2$ is bijective,
and we may choose an element  $b\in\maxid_B$ mapping to $\bar{b}$. 
Since the generic formal fiber of $B$ is reduced, the formal completion $\hat{B}$ is reduced,
by Proposition \ref{reduced, irreducible, integral}.
According to Corollary \ref{image canonical map}, there is an element
$a'\in \maxid_{\hat{R}}$ mapping to $\bar{b}$. 
Its class in the cotangent space $\bar{a}\in \maxid_{\hat{R}}/\maxid^2_{\hat{R}}$ is the image of some
element $a\in \maxid_R$.
Whence the map \eqref{cotangent map} on cotangent spaces is surjective.

Now choose some $a_1,\ldots, a_r\in\maxid_R$ so that their images in the cotangent
space form a vector space basis. According to the Nakayama Lemma, 
the elements $a_1,\ldots,a_r\in \maxid_B$ generate this maximal ideal.
It follows that the extension of local rings $R\subset B$ has invariant $e=1$.

According to Proposition  \ref{intermediate not finite}, 
the map $\hat{R}\ra\hat{B}$ is bijective and   $R\subset B$ is faithfully flat.
Furthermore, if the inclusion $R\subset B$ is not an equality, then $B$ is not a finite $R$-algebra.
\qed

\begin{corollary}
Suppose  the intermediate local ring $R\subset B\subset R^\prc$ is noetherian, and its generic
formal fiber is reduced.
If the local ring  $R$ is regular, or Cohen--Macaulay, the same holds for $B$.
\end{corollary}

\proof
A local noetherian ring is regular or Cohen--Macaulay
if and only if the respective property holds for its formal completion
(\cite{EGA IVa},  Proposition  17.1.5 and Proposition 16.5.2).
If $R$ is regular, or Cohen--Macaulay, the same holds for $\hat{R}=\hat{B}$, and thus for $B$.
\qed

\medskip
The assumption that the local ring $B$ is noetherian, with reduced generic formal fiber,
can be rephrased as a flatness condition:

\begin{proposition}
The intermediate ring $R\subset B\subset R^\prc$ is noetherian, with reduced generic formal fiber
if and only if the the extension of local rings  $B\subset \hat{R}$ is flat, and 
the induced map $\hat{B}\ra\hat{R}$ is injective.
\end{proposition}

\proof
The condition is necessary: The composite map
$$
\hat{R}\lra\hat{B}\lra\hat{R}
$$
is the identity,   the  map $\hat{R}\ra\hat{B}$ on the left is bijective by Theorem \ref{bijective}.
Whence the map $\hat{B}\ra\hat{R}$ on the right is even bijective.
Since $B\ra\hat{B}$ is flat, so must be the composition $B\ra\hat{B}\ra\hat{R}$.

The condition is also sufficient: The composite mapping 
$$
\Spec(\hat{R})\lra\Spec(B)\lra\Spec(R)
$$
is surjective, the map $\Spec(B)\ra\Spec(R)$ on the right is bijective according to Proposition \ref{universal homeomorphism},
whence the map $\Spec(\hat{R})\ra\Spec(B)$ on the left is surjective.
Therefore, the inclusion $B\subset\hat{R}$ is faithfully flat.
Now let $\idealb_0\subset\idealb_1\subset\ldots$ be an ascending chain of ideals in $B$.
The ascending chain $\idealb_0\hat{R}\subset\idealb_1\hat{R}\subset\ldots$ of ideals
in the noetherian ring $\hat{R}$
is stationary. 
By flatness, the canonical map $\idealb_i\otimes_B\hat{R}\ra\idealb_i\hat{R}$ is bijective, and 
by faithfully flat descent, the original ascending chain is stationary as well.
Thus the local ring $B$ is noetherian. Since $\hat{B}\ra \hat{R}$ is injective and $\hat{R}$ is reduced,
the ring $\hat{B}$ is reduced. By Proposition \ref{reduced, irreducible, integral}, the integral local ring $B$ has 
reduced generic formal fiber.
\qed

\medskip
Suppose from now on that the complete local noetherian ring $\hat{R}$ is normal.
According to Proposition \ref{intersection and closure}, also $R^\prc$ is normal.
Write $A_\lambda\subset R_\lambda\subset R^\prc$ for the normalization of the finite $R$-algebras $A_\lambda$.
This gives another filtered direct system $R_\lambda\subset R^\prc$, $\lambda\in L$
of subrings containing $R$, with $R^\prc=\bigcup_{\lambda\in L} R_\lambda$.
By the Mori--Nagata Theorem  
the $R_\lambda$ are Krull domains (see \cite{AC 7}, Chapter VII, \S1, No.\ 8, Proposition 12). 
Such rings are not necessarily noetherian.
Note, however, that if $R$ is one-dimensional, all the $R_\lambda$ are discrete valuation rings,
by the Krull--Akizuki Theorem (see \cite{AC 7}, \S2, No.\ 5, Corollary 2 of Proposition 5). If $R$ is two-dimensional,
the local rings $R_\lambda$ at least remain noetherian, by a result of Nagata \cite{Nagata 1955}
attributed to Mori.

We now want to give a sufficient condition for the $p$-radical closure $R^\prc$ to be
noetherian.  To proceed, consider the subset $L'\subset L$ of all indices $\lambda$ so that
the local Krull domain $R_\lambda $ is noetherian, and its generic formal fiber  
is  reduced. The latter happens, for example, if the $R_\lambda$ are regular.

\begin{theorem}
\mylabel{krull domain}
Suppose that all formal fibers of the local noetherian ring $R$ are reduced,
and that the subset $L'\subset L$ is cofinal.
Then the local ring $R^\prc$ is noetherian,
its generic formal fiber is reduced, the extension of local rings $R\subset R^\prc$
is faithfully flat, with  invariants $e=f=1$, and the canonical map $\hat{R}\ra\widehat{R^\prc}$
is bijective.
\end{theorem}

\proof
First note that the order set $L'$ is filtered,
because it is cofinal in the filtered ordered set $L$.
We thus have a filtered direct system $R_\lambda$, $\lambda\in L'$ of local noetherian
rings. All transition maps are faithfully flat, and local with invariants $e=f=1$.
Clearly,  $R^\prc=\bigcup_{\lambda\in L'} R_\lambda$.
It follows that the extension of local rings $R\subset R^\prc$ also has invariants $e=f=1$.
Consequently, the maximal ideal $\maxid_{R^\prc}=\maxid_RR^\prc$ is finitely generated.

The main issue is to verify that the ring $R^\prc$ is noetherian.
This means that every ideal is finitely generated.
According to Cohen's Theorem \cite{Cohen 1946}, it suffices to check this merely for the prime ideals  $\primid\subset R^\prc$.
Set $\primid_\lambda=\primid\cap R_\lambda$. For each $\lambda\leq \mu$, we have an inclusion
$\primid_\lambda R_\mu\subset\primid_\mu$, and  $\primid=\invlim \primid_\lambda$.
Since   $R_0=R$ is noetherian, it suffices to check that the inclusion 
$\primid_0 R_\mu\subset \primid_\mu$ is an equality
for all $ \mu\in L'$.
By Proposition \ref{universal homeomorphism}, the corresponding  closed embedding $\Spec(R_\mu/\primid_\mu)\subset\Spec (R_\mu/\primid_0 R_\mu)$
is a bijection.
It remains to verify that the scheme  $\Spec (R_\mu/\primid_0 R_\mu)$ is reduced.

Set $Y=\Spec(R_0)$ and $X=\Spec(R_\mu)$, write $f:X\ra Y$ for the canonical morphism,
with induced map $\hat{f}:\hat{X}\ra\hat{Y}$ on formal completions,
and consider the resulting commutative diagram of schemes
$$
\begin{CD}
\hat{Y} 	@<\hat{f} << 	\hat{X}\\
@Vq VV			@VVp V\\
Y		@<<f<	X.
\end{CD}
$$
According to Theorem \ref{bijective}, the upper vertical arrow is an isomorphism.
Let $C\subset Y$ be the integral closed subscheme defined by the ideal $\primid_0\subset R_0=R$.
Since the formal fibers of $Y$ are reduced, the scheme $q^{-1}(C)$ is reduced, whence 
also 
$$
p^{-1}(f^{-1}(C))=\hat{f}^{-1}(q^{-1}(C)).
$$
Since the morphism $p:\hat{X}\ra X$ is faithfully flat, the scheme $f^{-1}(C)=\Spec (R_\mu/\primid_0 R_\mu)$ must be reduced.

We next check that the extension of local rings $R\subset R^\prc$ has numerical invariants $e=f=1$.
The residual degree is $f=1$, which holds for all intermediate rings by Proposition \ref{universal homeomorphism}.
Since each extension of local rings $R\subset R_\lambda$, $\lambda\in L'$ has invariants $e=1$,
the same holds for the union $R^\prc=\bigcup_{\lambda\in L'} R_\lambda$.

Finally, we examine the generic formal fiber.
By assumption, the ring $R$ is formally integral. Proposition \ref{universal homeomorphism} 
ensures that $\dim(R)=\dim(R^\prc)$.
We thus may apply Proposition \ref{intermediate not finite} and infer that $\hat{R}\ra\widehat{R^\prc}$ is bijective.
In particular, the local noetherian ring $R^\prc$ is formally reduced.
Thus its generic formal fiber is reduced, according to Proposition \ref{reduced, irreducible, integral}.
\qed

\medskip
It is perhaps worthwhile to point out that there is a canonical exhaustive ascending chain
of intermediate rings $$R=B_0\subset B_1\subset \ldots \subset R^\prc$$  given by
\emph{heights}: Recall that  $F_\infty $ denotes  the relative $p$-radical closure
of $F=\Frac(R)\subset \Frac(\hat{R})$.
Let $F_n\subset F_\infty$ be the set of elements $a\in F_\infty$
of height $\leq n$, that is, of degree $[F(a):F]\leq p^n$,
and define $  B_n$ to be the integral closure of $R$ with respect to the field extension $F\subset F_n$.
Then we have $R^\prc=\bigcup_{n\geq 0} B_n$; it would be interesting to 
know whether these $B_n$ are noetherian, with reduced generic formal fiber.

\medskip 
In the   proof of Theorem \ref{reduction}, we have used the following basic fact:

\begin{lemma}
\mylabel{tensor product}
Let $B$ be a reduced $\FF_p$-algebra,   and $B\subset C$ be a  ring
extension such that every element $c\in C $ has $c^{p^\nu}\in B$ 
for some exponent $\nu\geq 0$. Then the homomorphism 
$$
C\lra (C\otimes_B C)_\red,\quad c\longmapsto (\text{\rm $c\otimes 1$ modulo nilradical})
$$
is bijective.
\end{lemma}

\proof
Let $f:C\ra C\otimes_B C$ be the homomorphism given by $c\mapsto c\otimes 1$,
and write $\bar{f}:C\ra (C\otimes_B C)_\red$ for the composite map.
The multiplication map $g:C\otimes_B C\ra C$ given by $c\otimes c'\mapsto cc'$
satisfies $g\circ f=\id_C$. It factors over $(C\otimes_B C)_\red$, because $C$ is reduced.
Thus $\bar{f}$ is injective.

To see that $\bar{f}$ is surjective, it suffices to show  that the class of each $1\otimes c$ 
lies in the image of $f$, up to nilpotent elements. Choose an exponent $\nu\geq 0$ with $c^{p^\nu}\in B$.
Using $\FF_p\subset B$, we have 
$ (c\otimes 1-1\otimes c)^{p^\nu}=0$. Thus $c\otimes 1-1\otimes c$ is nilpotent. From
$$
1\otimes c= f(c) - (c\otimes 1-1\otimes c) 
$$
we see that $\bar{f}$ is surjective.
\qed

%===========================================================
\section{Discrete valuation rings}
\mylabel{Discrete valuation rings}

We now specialize our results to the case that our local noetherian ring $R$ is a discrete valuation ring.
The field of fractions is given by  $F=\Frac(R)=R[1/a]$, where $a\in R$ is any non-zero non-unit.
Its formal completion $\hat{R}$ is an excellent discrete valuation ring. Obviously, there
are only two formal fibers: The closed formal fiber, which is the spectrum of the residue field 
$k=\hat{R}/\maxid_{\hat{R}}= R/\maxid_R$ and of little interest,
and the generic formal fiber, which is 
\begin{equation}
\label{ff field}
\hat{R}\otimes_R F=\hat{R}[1/a]=\Frac(\hat{R}).
\end{equation}
From this one deduces the following well-known fact:

\begin{proposition}
\mylabel{excellent dvr}
Consider the following three conditions:
\begin{enumerate}
\item The inclusion $R\subset R^\prc$ is bijective.
\item The field extension $F\subset\Frac(\hat{R})$ is separable.
\item The integral ring $R$ is japanese.
\item The discrete valuation ring $R$ is excellent.
\end{enumerate}
Then the implications (i) $\Leftarrow$ (ii) $\Leftrightarrow$ (iii) $\Leftrightarrow$ (iv) hold.
The four conditions are equivalent provided that  
the field extension $F\subset\Frac(\hat{R})$
can be written as  a purely inseparable extension followed by a separable extension.
\end{proposition}

\proof
First note that if $b\in F$ not contained in $R$, then $a=1/b\in R$ is a non-zero non-unit,
hence $R[b]=R[1/a]=F$. It follows that $R'=R$ and $R'=F$ are the only $R$-subalgebras $R'\subset F$.
In light of this, the implications (i) $\Leftarrow$ (ii) $\Leftrightarrow$ (iii) 
are special cases of Proposition \ref{geometrically reduced}.
The implication (iv) $\Rightarrow$ (ii) is trivial.
Conversely, if $F\subset\Frac(R)$ is separable, then   the two formal fibers of $R$ are geometrically regular.
The other two conditions for excellence also hold: The ring is obviously noetherian,
and since it is one-dimensional, it must be universally catenary, compare
\cite{EGA IVb}, Remark 7.1.2. Thus (ii) $\Rightarrow$ (iv) holds.
Under the additional assumption on the field extension $F\subset \Frac(\hat{R})$, the equivalence of all
four conditions again follows from Proposition \ref{geometrically reduced}.
\qed

\medskip
Let us examine the $p$-radical closure $R\subset R^\prc\subset\hat{R}$ of our discrete valuation ring
in more detail.
Recall that $R\subset A_\lambda\subset R^\prc$, $\lambda\in L$ denotes the filtered direct
system of all finite $R$-subalgebras. Each $A_\lambda$ is a one-dimensional integral local noetherian ring,
and the   extensions of local rings $R\subset A_\lambda$ is finite and flat.
We write $n_\lambda=\rank_R(A_\lambda)$ for the degrees.
 
\begin{proposition}
\mylabel{order extensions dvr}
The local extension $R\subset A_\lambda$ has numerical invariants
$f=1$ and $e=n_\lambda$. If  $n_\lambda\neq 1$, then the rings $R_\lambda$ are not 
normal.
\end{proposition}

\proof
The first assertion follows from Proposition \ref{universal homeomorphism}. Since $A_\lambda$ is finite and flat,
the formula $n=ef$ holds, and thus the ramification index coincides with the degree.
The last assertion follows from 
Corollary \ref{infinite normalization}.
\qed

\medskip
Let $A_\lambda\subset R_\lambda$ be the normalization inside the field $\Frac(\hat{R})$.
According to Proposition \ref{intersection and closure}, we have $R_\lambda\subset R^\prc$.
This gives another filtered direct system of $R$-subalgebras $R_\lambda\subset R^\prc$,
with $R^\prc=\bigcup_{\lambda\in L} R_\lambda$. Clearly, we have 
$$
n_\lambda=\rank_R(A_\lambda)=[R_\lambda\otimes_RF:F],
$$
where $F=\Frac(R)$ is the field of fractions.

\begin{proposition}
\mylabel{order extensions}
Each $R_\lambda$ is a discrete valuation ring, and the extension of local rings
$R\subset R_\lambda$  has   numerical invariants $e=f=1$ and $n=n_\lambda$.
If $n_\lambda\neq 1$, then $R_\lambda$ is not a finite $R$-algebra.
\end{proposition}

\proof
By the Krull--Akizuki Theorem (see \cite{AC 7}, \S2, No.\ 5, Corollary 2 of Proposition 5), 
the normalization $R_\lambda$ is a discrete valuation ring,
an in particular noetherian, with regular formal fibers. 
The other assertions thus follow from Theorem \ref{bijective}.
\qed

\begin{theorem}
Suppose that the field extension $F\subset\Frac(\hat{R})$ can be written as a
purely inseparable extension followed by a separable extension.
Then the $p$-radical closure  $R^\prc$ is an excellent discrete valuation ring. 
The extension of local rings $R\subset R^\prc$ has invariants $e=f=1$ and $n=\sup_{\lambda\in L}\{n_\lambda\}$,
and the induced map on formal completions $\hat{R}\ra\hat{R}^\prc$ is bijective.
If $R$ is not excellent, then the  ring extension $R\subset R^\prc $ is integral but not finite.
\end{theorem}

\proof
Obviously, for every index $\lambda\in L$, the local Krull domain $A_\lambda'=R_\lambda$ is noetherian and regular,
in particular its formal fibers are reduced.
Thus Theorem \ref{krull domain} applies, and the result follows.
\qed

\medskip
We now can give the following universal property:

\begin{corollary}
Assumptions as in the theorem.
Let $R\subset A$ be an extension of local rings.
If $A$ is an excellent discrete valuation ring, or more generally any formally integral  local noetherian ring
with geometrically reduced generic formal fiber,
then we have $R^\prc\subset A$ inside the field
of fractions $\Frac(\hat{A})$.
\end{corollary}

\proof
Since $R$ is a discrete valuation ring, the induced map $\hat{R}\ra\hat{A}$ is injective:
Otherwise it would factor over $\hat{R}/\maxid_{\hat{R}}$, and the image of 
$\Spec(\hat{A})\ra\Spec(\hat{R})$ consists of the closed points.
The same then holds for composite map $\Spec(\hat{A})\ra\Spec(A)\ra\Spec(R)$.
The latter maps, however, are dominant, contradiction.
We thus may apply Corollary  \ref{minimal}, and the assertion follows.
\qed

\medskip
Let me close this paper with the following question:
Do similar results as for discrete valuation rings hold for 
regular local noetherian rings?

%===========================================================

\end{document}